\documentclass[AMA,STIX1COL]{WileyNJD-v2}
\usepackage{amsmath}

\newcommand{\der}[2]{\frac{\partial{#1}}{\partial{#2}}}
\newcommand{\dder}[2]{\frac{\partial^2{#1}}{\partial{#2}^2}}
\newcommand{\norm}[1]{\left\lVert#1\right\rVert}
\newcommand{\erfc}{\text{erfc}}
\newtheorem{exm}{Example}[section]

\articletype{Research Paper}%

\received{}
\revised{}
\accepted{}

\begin{document}

\title{Integral Transform Solution of Random Coupled Parabolic Partial Differential Models}

\author[1]{Mar\'{i}a Consuelo	Casab\'{a}n*}

\author[1]{Rafael Company}

\author[2]{Vera N. Egorova}

\author[1]{Lucas J\'{o}dar}
 
\authormark{ M.-C. Casab\'{a}n \textsc{et al}}

\address[1]{\orgdiv{Instituto de Matem\'{a}tica Multidisciplinar}, \orgname{Universitat Polit\`{e}cnica de Val\`{e}ncia}, \orgaddress{\state{Camino de Vera, s/n, 46022 Valencia}, \country{ Spain}}}

\address[2]{\orgdiv{Depto. de Matem\'{a}tica Aplicada y Ciencias de la Computaci\'{o}n}, \orgname{Universidad de Cantabria}, \orgaddress{\state{Avda. de los Castros, s/n, 39005 Santander}, \country{Spain}}}

\corres{*Mar\'{i}a Consuelo	Casab\'{a}n*. \email{macabar@imm.upv.es}}


\abstract[Summary]{	
	Random coupled parabolic partial differential models are solved numerically using random cosine Fourier transform together with non Gaussian random numerical integration that capture the highly oscillatory behavior of the involved integrands. Sufficient  condition of spectral type imposed on the random matrices of the system are given so that the approximated stochastic process solution and its statistical moments are numerically convergent. Numerical experiments illustrate the results.
}



\maketitle
%

\section{Introduction}

Random time dependent scalar mean square partial differential models have been treated recently from both the theoretical and numerical points of view \cite{}, because  in real problems the parameters, coefficients and initial/boundary conditions are subject to uncertainties, not only by error measurement but also due to the difficulty of access to the measurement, the possible heterogeneity of the materials or media, etc.
Spatial uncertainty models described by random elliptic PDEs in bounded domains are treated in \cite{Back2011, Bachmayr2017, Ernst2018} using spectral Galerkin and collocation methods.
 Dealing the coupled partial differential models the uncertainties are involved in the matrix coefficients or vector initial/boundary conditions.

Coupled partial differential models are frequent in several engineering disciplines such as geomechanics \cite{Sheng_1995},  geotechnics \cite{Mitchell_1991},
 microwave heating processes \cite{Metaxas_1983}, optics \cite{Das_1991}, ocean models \cite{Ekman_1905,Yosef2017}, etc. They also appear in plasma fusion models \cite{Weston_1981}, cardiology \cite{Hodgkin_1952,Winfree_1987}, or species population dynamics \cite{Galiano2012}.
 
 Solving random models presents somewhat unexpected peculiarities not presented in the deterministic case. In fact, in the random case, it is important not only the determination of the exact or approximate stochastic process solution but also the computation of its statistical moments, mainly the expectation and the standard deviation.
 
 Using iterative methods involve the storage of high computational complexity of all the previous levels of iteration and usually the methods become unmanageable \cite{Casaban2019b}. 
 This motivates the search of alternative methods with as simple as possible expressions of the stochastic solution process. 
 In the recent paper \cite{Casaban2019a}, one uses Fourier transforms together with the random Gaussian quadrature rules to approximate stochastic process solution. The method proposed in \cite{Casaban2019a} has the advantage that approximate stochastic solution is simple and the computation of its statistical moments are manageable, however, as the numerical integration is based on Gaussian quadrature rules, the accuracy decreases for highly oscillatory Fourier kernels in large domains, \cite{Shampine2008,Iserles2004, Ma2018}. In this paper, the above mentioned drawback is overcome by taking an appropriate truncation of the infinite integral and using quadrature rules with good behaviour dealing with highly oscillatory integrands.
 
 PDE models in unbounded domains using Fourier integral transforms have been treated in \cite{JodarGoberna_1996,JodarGoberna_1998}. This paper deals with a more general random coupled parabolic problem
 \begin{align}
 &\der{u(z,t)}{t}(\xi) = A(\xi) u(z,t)(\xi) + B(\xi) \dder{u(z,t)}{z}(\xi), \quad z>0, \quad t>0\,,\label{eq:problem1}\\ 
 &u(z,0)(\xi) = f(z)(\xi), \quad z>0, \; \xi \in \Omega,
 \\
 &\der{u}{z} (0,t)(\xi) = g(t)(\xi), 
 \quad  t>0, \; \xi \in \Omega, \label{eq:problem3}\\
& \lim_{z\rightarrow \infty} u(z,t)(\xi) = 0, \quad \lim_{z\rightarrow \infty} \der{u}{z}(z,t)(\xi) = 0, \label{eq:problem4}
 \end{align}
 where $u(z,t)(\xi) = [u_1(z,t)(\xi), u_2(z,t)(\xi)]^T \in \mathbb{R}^2$,  $g(t)(\xi) = [g_1(t)(\xi), g_2(t)(\xi)]^T $, $f(z)(\xi) = [f_1(z)(\xi), f_2(z)(\xi)]^T $, and 
\begin{equation}
 A(\xi) = \left(a_{ij}(\xi)\right)_{1 \leq i,j \leq 2}, \; B(\xi) = \left(b_{ij}(\xi)\right)_{1 \leq i,j \leq 2}.
\end{equation}

Here $A(\xi)$ and $B(\xi)$ are random matrices, $f(z)(\xi)$ and  $g(t)(\xi)$ are stochastic processes (s.p.'s) with properties to be specified later. We assume that for each event $\xi \in \Omega$, the sample matrix $B(\xi)$ satisfies
 
 \begin{equation}\label{eq:lambda_min}
 \lambda_{\min} \left( \frac{B(\xi)+B(\xi)^T}{2} \right) = b(\xi) > 0,
 \end{equation}
 where $ \lambda_{\min}$ denotes the minimum eigenvalue. 
 
 This paper is organized as follows. Section \ref{sec:determ} deals with the solution of the simplified deterministic problem after taking sample realizations for each event $\xi \in \Omega$. Matrix analysis of involved matrices $A(\xi)$ and $B(\xi)$ is performed in order to determine the spectral sufficient condition to guarantee the convergence. 
 
 The unsuitable use of Gaussian quadrature for cosine oscillatory integrals, and the convergence of the truncated integrals suggest the introduction of alternative quadrature formulae such as the midpoint Riemann sum, see \cite{Davis&Rabinowitz} section 3.9. 
 
 In Section \ref{sec:random}, the random case is addressed taking into account the ideas of previous section in order to construct random approximate solution s.p.'s that make manageable the computation of its statistical moments, in particular, the expectation and the variance. Simulations show the efficiency of the proposed numerical methods.

 \section{Solving the Sampled Deterministic Case}\label{sec:determ}

For the sake of clarity in the presentation, let us recall some algebraic concepts, notations and results. 

If $P$ is a matrix in $\mathbb{R}^{N \times N}$, its logarithmic operator norm $\mu(P)$ is defined by
\begin{equation}
\mu(P) = \max \left\lbrace \lambda; \lambda \text{ eigenvalue of } \frac{P + P^T}{2} \right\rbrace.
\end{equation}

By \cite{Dahlquist1958}, the matrix exponential $e^{Pt}$ satisfies 
$
\norm{e^{Pt}} \leq e^{t \mu(P)}, \quad t \geq 0.
$

\begin{lemma}\label{lemma:1}
	Let $B \in \mathbb{R}^{N \times N}$ be a matrix such that $B + B^T$ is positive definite and satisfies (\ref{eq:lambda_min}). Then
	
	\begin{equation} 
	\mu(A - \omega^2 B) \leq \mu(A) - b\omega^2, \quad b = \lambda_{\min} \left( \frac{B+B^T}{2} \right), \quad \omega>0\,.
	\end{equation}
\end{lemma}

\begin{proof}
	Let $\omega >0$, and let us write
\begin{equation} \label{eq:matrixAAtBBtw}
\frac{(A-B\omega^2)+ (A-B\omega^2)^T}{2} = \frac{A + A^T}{2} - \omega^2 \left( \frac{B+B^T}{2} \right).
\end{equation}

	As $\frac{A + A^T}{2}$ and $ - \omega^2 \left( \frac{B+B^T}{2} \right)$ are both symmetric matrices, by Ostrowski theorem, see \cite{Ostrowski1959}, each eigenvalue $\lambda$ of matrix \eqref{eq:matrixAAtBBtw}  satisfies
	
	\begin{equation}
	\lambda \leq \lambda_{\max}\left( \frac{A+A^T}{2} \right) - \omega^2 \lambda_{\min}\left( \frac{B+B^T}{2} \right) = \mu(A) - b \omega^2, 
	\end{equation}	
	where $b$ is given by (\ref{eq:lambda_min}). Hence, the result is established.
\end{proof}

If $f(z) = [f_1(z), f_2(z)]^T$ is absolutely integrable in $[0, \infty)$, then the cosine Fourier transform of $f(z)$ is defined by
\begin{equation}
\mathcal{F}_c[f](\omega) = \int_{0}^{\infty} f(z) \cos (\omega z) dz, \quad \omega \geq 0,
\end{equation}
and if $f(z)$ is twice differentiable and $f''(z) = [f_1''(z), f_2''(z)]^T$ is absolutely integrable, then

\begin{equation} \label{eq:fourier_der}
\mathcal{F}_c[f''](\omega) = -\omega^2 \mathcal{F}_c[f](\omega) - f'(0), \quad \omega \geq 0.
\end{equation}

In order to simplify the notation in the following content of this section, we will denote for a realization of random matrices $A(\xi)$ and $B(\xi)$ as $A$ and $B$, respectively. 
For the sake of coherence we also denote $f(z)$, $g(t)$ and $u(z,t)$ as the corresponding realizations for a fixed event $\xi \in \Omega$.

In order to obtain a candidate solution of problem (\ref{eq:problem1})--(\ref{eq:lambda_min}), let us apply the cosine Fourier transform
$\mathcal{F}_c$ regarding $u = u(\cdot, t)$ as an absolute integrable function of the active variable $z>0$. Let us denote 

\begin{equation}
V(t) = \mathcal{F}_c[u(\cdot,t)](\omega) = \int_{0}^{\infty} u(z,t) \cos(\omega z) dz, \quad \omega >0,
\end{equation}
and applying $\mathcal{F}_c$ to (\ref{eq:problem1}) and taking into account (\ref{eq:fourier_der}), one gets

\begin{align}
&\mathcal{F}_c\left[\dder{u}{z}(\cdot,t)\right](\omega) = -\omega^2 \mathcal{F}_c[u(\cdot,t)] (\omega)- \der{u}{z}(0,t) = -\omega^2 V(t)(\omega)-g(t),
\\
&V(0) = \mathcal{F}_c[u(0,t)](\omega)  = \mathcal{F}_c[
f(z)](\omega) = F(\omega).
\end{align}

Hence, $V(t)(\omega)$ is the solution of the initial value problem in time,

\begin{equation}
\label{eq:ode}
\frac{d}{dt} V (t) (\omega) = (A-\omega^2 B)V(t) (\omega) - B g(t), \quad t>0, \; \omega>0 \text{ fixed}, \quad V(0)(\omega) = F(\omega).
\end{equation}

The solution of (\ref{eq:ode}) for $\omega>0$ fixed, takes the form

\begin{equation}  \label{eq:cs}
V(t)(\omega) = e^{(A-\omega^2 B)t} \left\lbrace F(\omega) - \int_{0}^{t} e^{-(A-\omega^2B)s} c(s)ds \right\rbrace, \quad  c(s) = Bg(s).
\end{equation}

Under hypothesis of Lemma \ref{lemma:1}, and continuity of $g(t)$, taking cosine inverse $\mathcal{F}_c^{-1}$, one gets

\begin{equation}\label{eq:u_def}
u(z,t) = \mathcal{F}_c^{-1}[V(t)(\omega)] = \frac{2}{\pi} \int_{0}^{\infty} e^{(A-\omega^2B)t} F(\omega) \cos(\omega z) d \omega = \frac{2}{\pi} (I_1 - I_2),
\end{equation}
where 

\begin{equation}\label{eq:i1}
I_1 = \int_{0}^{\infty} e^{(A-\omega^2B)t} F(\omega) \cos(\omega z) d\omega,
\quad I_2 = \int_{0}^{\infty} \left( \int_{0}^{t}e^{(A-\omega^2B)(t-s)} c(s) \cos(\omega z) ds \right)d\omega.
\end{equation}

Integrals  (\ref{eq:i1}) can be truncated for $\omega$, getting the approximations

\begin{equation}\label{eq:i1r}
I_1(R) = \int_{0}^{R} e^{(A-\omega^2B)t} F(\omega) \cos(\omega z) d\omega, 
\quad
I_2(R) = \int_{0}^{R} \left( \int_{0}^{t}e^{(A-\omega^2B)(t-s)} c(s) \cos(\omega z) ds \right)d\omega,  \quad R>0,
\end{equation}
and the approximate solution $
u_R (z,t) = \frac{2}{\pi} \left(I_1(R) - I_2(R)\right), \; z>0,\; t>0.
$

Now we prove that $\left\lbrace u_R(z,t)\right \rbrace$ is convergent and that
$
\lim_{R\rightarrow \infty} u_R(z,t)= u(z,t),
$ stating that

\begin{equation}\label{eq:u_lim}
u(z,t)- u_R(z,t) = \frac{2}{\pi} \left(J_1(R) - J_2(R)\right) \xrightarrow[R \rightarrow \infty]{} 0,
\end{equation}
where

\begin{equation}\label{eq:j1r}
J_1(R) = \int_{R}^{\infty} e^{(A-\omega^2B)t} F(\omega) \cos(\omega z) d\omega, 
\quad
J_2(R) = \int_{0}^{t} \left( \int_{R}^{\infty}e^{(A-\omega^2B)(t-s)} c(s) \cos(\omega z) d\omega \right)ds.
\end{equation}

Using Lemma \ref{lemma:1}, (\ref{eq:j1r}) and the substitution 
$u=\omega\sqrt{bt}$ one gets

\begin{equation}\label{eq:j1_bound}
\begin{split}
\norm{J_1(R)} &\leq \int_{R}^{\infty} e^{\mu (A-\omega^2B)t} \norm{F(\omega)} d\omega \leq 
\norm{F}_{\infty}  \int_{R}^{\infty} e^{(\mu (A)-\omega^2b)t}d\omega
= \norm{F}_{\infty} e^{\mu (A)t}  \int_{R}^{\infty} e^{-bt\omega^2}d\omega 
\\
&= \frac{  \norm{F}_{\infty} e^{\mu (A)t}}{\sqrt{bt}}\int_{R\sqrt{bt}}^{\infty} e^{-u^2}du = \frac{  \norm{F}_{\infty} \sqrt{\pi} }{2 \sqrt{bt}}e^{\mu (A)t} \erfc(R\sqrt{bt}),
\end{split}
\end{equation}
where $
\sup \left\lbrace \norm{F(\omega)}; \; \omega\geq 0\right\rbrace =  \norm{F}_{\infty}.
$

Also, from Lemma \ref{lemma:1}, for the second integral of  (\ref{eq:j1r}) one gets

\begin{equation}\label{eq:j2_prel1}
\norm{J_2(R)} \leq \int_{0}^{t} \left( \int_{R}^{\infty}e^{(\mu(A)-\omega^2b)(t-s)} \norm{B} \norm{g(s)} d\omega \right)ds
= \norm{B} \int_{0}^{t} \norm{g(s)}e^{\mu(A)(t-s)} \left( \int_{R}^{\infty}e^{-\omega^2b(t-s)} d\omega \right)ds.
\end{equation}

As we did in previous bound of $J_1(R)$ in (\ref{eq:j1_bound}), we have

\begin{equation}\label{eq:j2_prel2}
\int_{R}^{\infty} e^{-\omega^2 b (t-s)}d\omega = \frac{  \sqrt{\pi} }{2\sqrt{b(t-s)}}\erfc\left(R\sqrt{b(t-s)}\right). \end{equation}

From (\ref{eq:j2_prel1}) and (\ref{eq:j2_prel2}) one concludes

\begin{equation}\label{eq:j2_bound}
\norm{J_2(R)} \leq \frac{  \sqrt{\pi} \norm{B} }{2\sqrt{b}}
\int_{0}^{t} \frac{\norm{g(s)}e^{\mu(A)(t-s)}}{\sqrt{t-s}}  \erfc\left(R\sqrt{b(t-s)}\right) ds  = \frac{  \sqrt{\pi} \norm{B} }{b}
\int_{0}^{\sqrt{bt}} \norm{g\left(t- \frac{v^2}{b}\right)}e^{\mu(A)\frac{v^2}{b}}  \erfc\left(Rv\right) dv.
\end{equation}

As $\lim_{x\rightarrow \infty} \erfc(x) = 0$, and $g(t)$ is continuous and bounded in a bounded interval, from (\ref{eq:j1_bound}) and (\ref{eq:j2_bound}) it follows that

\begin{equation}\label{eq:j_lim}
\lim_{R\rightarrow \infty} J_i(R) = 0, \quad i=1,2.
\end{equation}

Hence, from (\ref{eq:u_lim}) and (\ref{eq:j_lim})
one gets $
\lim_{R\rightarrow \infty}( u(z,t)- u_R(z,t)) = 0.
$

Note also that $u(z,t)$ given by (\ref{eq:u_def})--(\ref{eq:i1}) is well defined because integrals $I_1$ and $I_2$ of (\ref{eq:i1}) are absolutely integrable, because for $z>0$ and $t>0$ fixed we have, see Lemma \ref{lemma:1},

\begin{equation*}
\norm{I_1} \leq \int_{0}^{\infty} e^{\mu(A)t} \norm{F(s)} e^{-\omega^2 bt} d \omega
\leq \norm{F}_{\infty} e^{\mu(A)t}  \int_{0}^{\infty}e^{-\omega^2 bt} d\omega < + \infty,
\end{equation*}
\begin{equation*}
\begin{split}
\norm{I_2} &\leq \int_{0}^{t} \int_{0}^{\infty} e^{\mu(A)(t-s)- \omega^2 b(t-s)} \norm{c(s)} d\omega ds 
\leq \norm{B} \int_{0}^{t} \norm{g(s)} e^{\mu(A)(t-s)} \left(\int_{0}^{\infty} e^{-\omega^2 b (t-s)} d\omega\right) ds
\\
&= \frac{\norm{B}}{\sqrt{b}} \int_{0}^{t} \frac{\norm{g(s)} e^{\mu(A)(t-s)}}{\sqrt{t-s}} \left(\int_{0}^{\infty} e^{-v^2}dv\right)ds
= \frac{\sqrt{\pi} \norm{B}}{2\sqrt{b}} \int_{0}^{t} \frac{\norm{g(s)}e^{\mu(A)(t-s)}}{\sqrt{t-s}}ds.
\end{split}
\end{equation*}

Summarizing the following result has been established:

\begin{theorem} \label{the:01}
	Consider the problem (\ref{eq:problem1})--(\ref{eq:problem4})  for a fixed event $\xi\in{\Omega}$, where $f(z)$ is absolutely integrable, $g(t)$ is continuous and $B$ satisfies  condition (\ref{eq:lambda_min}) with $b>0$. Then
	\begin{enumerate}[(i)]
		\item $u(z,t)$ given by (\ref{eq:u_def})--(\ref{eq:i1}) is solution of problem (\ref{eq:problem1})--(\ref{eq:problem4});
		\item $u_R (z,t) = \frac{2}{\pi} \left(I_1(R) - I_2(R)\right), \; z>0,\; t>0$, where $I_1(R)$ and $I_2(R)$ are defined by \eqref{eq:i1r}, converges as $R \rightarrow \infty$ to $u(z,t)$ uniformly for $z>0$ and point-wise at each $(z,t) \in \mathbb{R}^{+} \times \mathbb{R}^{+}$. 
	\end{enumerate}
\end{theorem}


\begin{exm}\label{exm:determ}
	Consider problem (\ref{eq:problem1})--(\ref{eq:problem4}) for a realization corresponding to $\xi_0 \in \Omega$ fixed with the data
	\begin{equation}
	A = \begin{bmatrix}
	0 & a \\-a & 0
	\end{bmatrix}, \quad 
	B = \nu I = \begin{bmatrix}
	\nu & 0\\ 0 & \nu
	\end{bmatrix},
	\quad a>0, \; \nu >0;
	\quad
	u(z,0) = f(z) = \begin{bmatrix}
	0\\0
	\end{bmatrix}, \quad \der{u}{z}(0,t) = \begin{bmatrix}
	-g(t)\\0
	\end{bmatrix},
	\end{equation}
	modelling the influence of the Earth's rotation on ocean currents \cite{Ekman_1905}, whose exact solution takes the form
	
	\begin{equation}\label{eq:exact}
	u_1(z,t) = \sqrt{\frac{\nu}{\pi}} \int_{0}^{t} \frac{g(s)}{\sqrt{t-s}} e^{-\left(\frac{z^2}{4\nu(t-s)}\right)} \cos \left(a(t-s)\right)ds,
	\quad u_2(z,t) = -\sqrt{\frac{\nu}{\pi}} \int_{0}^{t} \frac{g(s)}{\sqrt{t-s}} e^{-\left(\frac{z^2}{4\nu(t-s)}\right)} \sin \left(a(t-s)\right)ds.
	\end{equation}
\end{exm}
Here variables $t$ and $z$ represent the time and depth coordinates, $u_1$ and $u_2$ describe the zonal and meridional surface ocean current velocities, $a$ is the coriolis parameter and $\nu$ is the eddy parameterized vertical viscosity coefficient.

Note that in this case
\begin{equation}   \label{eq:AwB}
A-\omega^2 B
= \begin{bmatrix}
-\omega^2 \nu & a\\ -a & -\omega^2 \nu
\end{bmatrix}, 
\qquad \mu(A-\omega^2B) = \max\left\lbrace\lambda \in \sigma\begin{bmatrix}
-\omega^2 \nu & 0\\ 0 & -\omega^2 \nu
\end{bmatrix}\right\rbrace = -\omega^2 \nu.
\end{equation}

From \eqref{eq:cs} and \eqref{eq:AwB} it follows that 
\begin{equation}  \label{eq:cs-AwB}
c(s) = \nu
\begin{bmatrix}
-g(t)\\ 0
\end{bmatrix}; \quad
e^{(A-\omega^2 B)t} = e^{-\omega^2 \nu t}\begin{bmatrix}
\cos(at) & \sin(at)\\ -\sin(at) & \cos(at)
\end{bmatrix}.
\end{equation}

By (\ref{eq:u_def})--(\ref{eq:i1}) and \eqref{eq:cs-AwB} the solution of the problem takes the form

\begin{equation}\label{eq:int_solution}
u (z,t) = - \frac{2}{\pi} \int_{0}^{\infty} \int_0^t e^{(A-\omega^2 B)(t-s)} c(s) \cos(\omega z) ds d\omega = \frac{2}{\pi} \int_{0}^{\infty} I_2 (t, \omega) \cos (\omega z)d\omega,
\end{equation}
where 

\begin{equation}\label{eq:I2}
I_2(t, \omega) = \int_0^t e^{-\omega^2 \nu (t-s)}  \nu g(s) \begin{bmatrix}
\cos\left(a(t-s)\right)\\ -\sin\left(a(t-s)\right)
\end{bmatrix}ds.
\end{equation}

Taking advantage of the knowledge of the exact solution given by (\ref{eq:exact}) we check that Gauss-Laguerre quadrature \cite{Casaban2019a} of (\ref{eq:int_solution}) provides wrong results for large values of $z$.

Let us take the data $a = 1$, $\nu = 1$ and $g(t) = 1$. Note that in the case $g(t) =1$ expression (\ref{eq:I2}) becomes
\begin{equation}
I_2(t, \omega) = \frac{\nu}{a^2+\nu^2 \omega^4} \begin{bmatrix}
\omega^2 \nu + e^{-\omega^2 \nu t} \left( a \sin (at) - \omega^2 \nu \cos (at)\right)
\\
-a+ e^{-\omega^2 \nu t} \left( a \cos (at) + \omega^2 \nu \sin (at)\right)
\end{bmatrix}.
\end{equation}

Next Table \ref{table:gauss} shows the absolute errors when one approximates (\ref{eq:exact}) using Gauss-Laguerre quadrature of several degrees $M$ for $z= 5$ and $t= 1$.

\begin{table}
\centering
\begin{tabular}{|c|cc|}
\hline
$M$	&	AbsErr$\left(u_1(5,1)\right)$	&	 AbsErr$\left(u_2(5,1)\right)$\\
\hline
$ 1$	& $2.7254e-01$ & 	 $1.2057e-01$ \\ 	
$ 2$	& $7.4014e-01$ & 	 $3.5587e-01$ \\ 	
$ 3$	& $2.2645e-02$ & 	 $1.0762e-01$ \\ 	
$ 4$	& $4.4318e-01$ & 	 $9.7300e-02$ \\ 	
$ 5$	& $4.5831e-01$ & 	 $1.5709e-01$ \\ 	
$ 6$	& $3.6717e-01$ & 	 $2.1590e-01$ \\ 	
$ 7$	& $5.1360e-01$ & 	 $1.8173e-01$ \\ 	
$ 8$	& $1.9483e-03$ & 	 $5.3340e-02$ \\ 	
$ 9$	& $1.3911e-01$ & 	 $6.4812e-02$ \\ 	
$ 10$	& $4.8348e-01$ & 	 $1.3559e-01$ \\ 	
$ 11$	& $3.3161e-01$ & 	 $1.2510e-01$ \\ 	
$ 12$	& $2.1783e-01$ & 	 $6.6086e-02$ \\ 	
$ 13$	& $1.4817e-01$ & 	 $8.2631e-03$ \\ 	
$ 14$	& $2.3801e-01$ & 	 $5.7655e-02$ \\ 	
$ 15$	& $2.8347e-01$ & 	 $7.0344e-02$ \\ 	
\hline
\end{tabular}
\caption{Absolute error of numerical approximation of (\ref{eq:exact}) by using Gauss-Laguerre quadratures of degree $M$ at $z= 5$ and $t = 1$.}
\label{table:gauss}
\end{table}

The convergence of truncated integrals in Theorem \ref{the:01} suggests the approximation of the  truncated integrals using appropriate quadrature rules preserving the oscillatory behaviour. Let us denote
\begin{equation}\label{eq:uR}
u_R(z,t) = \frac{2}{\pi} \int_0^R I_2 (t, \omega) \cos(\omega z)d\omega.
\end{equation}

The midpoint Riemann sum proposed in \cite{Davis&Rabinowitz}, Section 3.9, applied to (\ref{eq:uR}) gives the approximation

\begin{equation}\label{eq:midpoint}
u_R(z,t) \approx \frac{2h}{\pi} \sum_{j= 0}^{N-1} I_2 \left(t, \left(j+\frac{1}{2}\right)h\right) \cos \left(\left(j+\frac{1}{2}\right)hz \right),
\end{equation} 
where $Nh = R$.

The numerical convergent of the proposed midpoint Riemann sum approximation of the solution given by \eqref{eq:midpoint} is studied with respect to the involved parameters $R$,  $N$ and $h$ in Table \ref{table:fixed_h} and Table \ref{table:fixed_r}.  In particular, in Table \ref{table:fixed_h} we fix the value of $h$ and vary $R$ in order to analyse the impact of the truncation point $R$. As expected, increasing values of $R$ result in smaller absolute error.  In Table \ref{table:fixed_r}, $R$ is fixed and the step-size discretization $h$
 is changing.

\begin{table}
	\centering
	\begin{tabular}{|c|cc|}
		\hline 
		$R$	& AbsErr$\left(u_1(5,1)\right)$ &  AbsErr$\left(u_2(5,1)\right)$\\
		\hline
$ 5$	& $2.2187e-04$ & 	 $6.0459e-06$ \\ 	
$ 10$	& $6.6113e-05$ & 	 $1.1396e-06$ \\ 	
$ 15$	& $9.6916e-05$ & 	 $4.9150e-07$ \\ 	
$ 20$	& $9.3665e-05$ & 	 $2.4768e-07$ \\ 	
$ 25$	& $8.4522e-05$ & 	 $1.3928e-07$ \\ 	
$ 30$	& $7.4942e-05$ & 	 $8.4709e-08$ \\ 		
		\hline
	\end{tabular}
	\caption{Absolute errors of numerical approximations of (\ref{eq:uR}) at $z=5$, $t = 1$, by the midpoint Riemann sum with fixed $h=0.05$. }
	\label{table:fixed_h}
\end{table}

\begin{table}
	\centering
	\begin{tabular}{|c|cc|}
		\hline 
		$h$	& AbsErr$\left(u_1(5,1)\right)$ &  AbsErr$\left(u_2(5,1)\right)$\\
		\hline
$ 0.2000$	& $1.4793e-04$ & 	 $3.4924e-07$  \\ 	
$ 0.1000$	& $1.4441e-05$ & 	 $5.0984e-08$  \\ 	
$ 0.0500$	& $9.3665e-05$ & 	 $2.4768e-07$  \\ 	
$ 0.0250$	& $1.3112e-04$ & 	 $3.4100e-07$  \\ 	
$ 0.0125$	& $1.4912e-04$ & 	 $3.8592e-07$  \\ 	
		\hline
	\end{tabular}
	\caption{Absolute errors of numerical approximations of (\ref{eq:uR}) at $z=5$, $t = 1$, by the midpoint Riemann sum with fixed $R = 20$ and various step-size $h$. }
	\label{table:fixed_r}
\end{table}

\begin{figure}[h]
	\begin{center}
			\centering
			\includegraphics[width=0.5\linewidth]{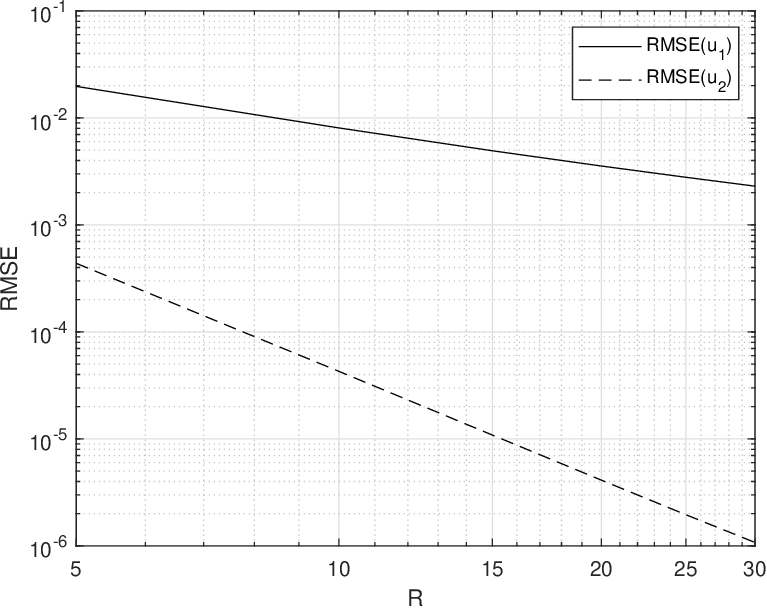}
			\caption{RMSE of numerical approximations of (\ref{eq:uR}) by the midpoint Riemann sum with fixed $h=0.05$ for various values of $R$ in the domain $(z,t)\in{[0,5]\times[0,1]}$ for the step-sizes $\Delta z=0.05$ and $\Delta t=0.01$.}
			\label{fig:test1_rmse}
	\end{center}
\end{figure}

\begin{figure}[h]
	\begin{center}
		\begin{minipage}[h]{0.45\linewidth}
			\centering
			\includegraphics[width=\linewidth]{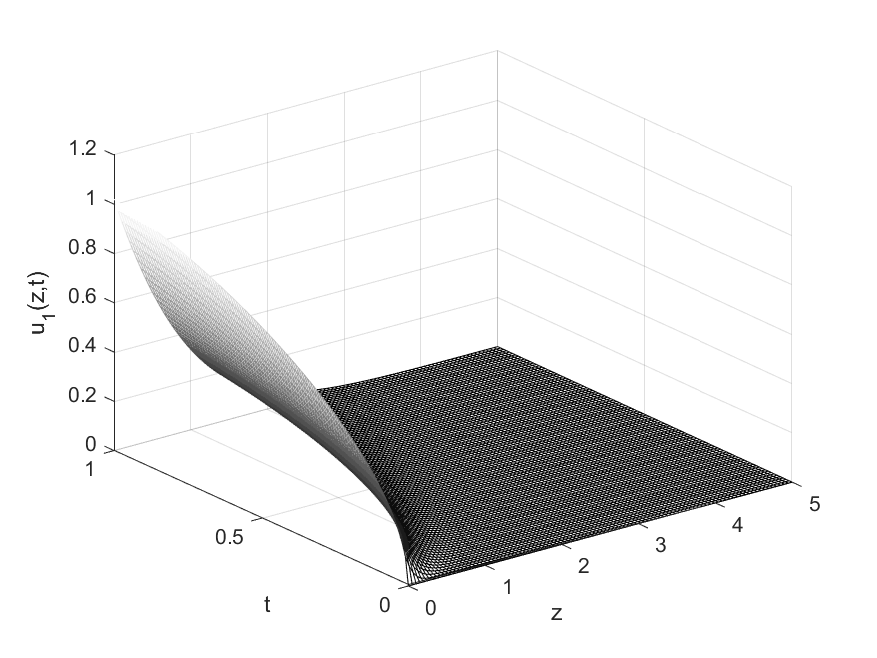}
			\caption{Solution  $u_1(z,t)$ calculated by (\ref{eq:exact}) for $a=\nu = 1, \; g(t)=1$.}
			\label{fig:test1_exact1}
		\end{minipage}
		\hfill
		\begin{minipage}[h]{0.45\linewidth}
			\centering
			\includegraphics[width=\linewidth]{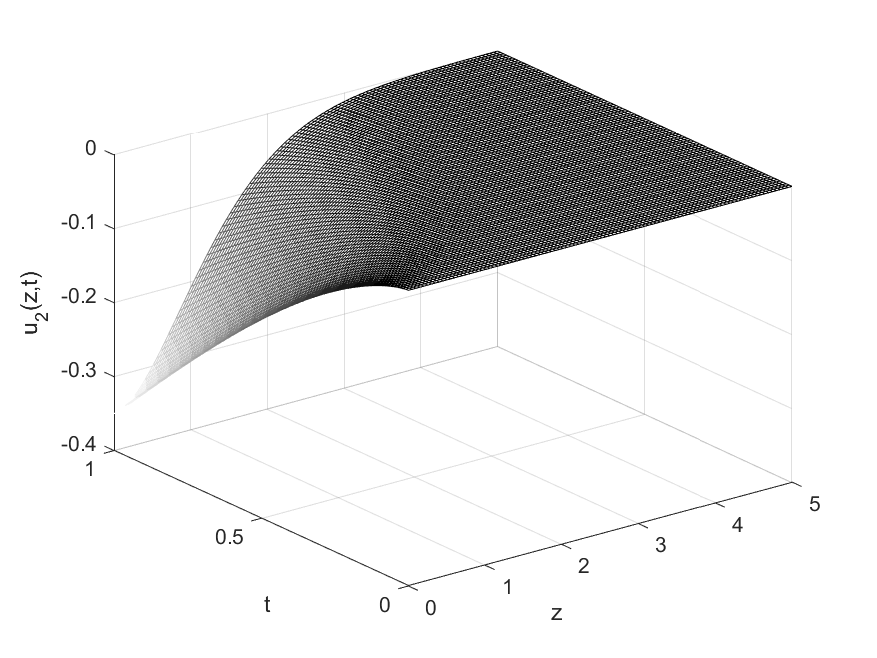}
			\caption{Solution  $u_2(z,t)$ calculated by (\ref{eq:exact}) for $a=\nu = 1, \; g(t)=1$.}
			\label{fig:test1_exact2}
		\end{minipage}
	\end{center}
\end{figure}

In order to compare the numerical approximation by the Riemann midpoint sum with the analytical solution (\ref{eq:exact}) in an entire computational domain, we calculate the root mean square error (RMSE). The results for fixed $h = 0.05$ at the fixed computational domain $0 \leq t \leq 1$, $0 \leq z \leq 5$ are presented in Table \ref{table:rmse_fixed_h} and plotted in Figure \ref{fig:test1_rmse}. Since the step-size is fixed, the total computational time varies depending on the size of the integration domain, i.e., on the values of $R$. The CPU time of each simulation is reported  in Table \ref{table:rmse_fixed_h} as well. The solutions $u_1(z,t)$  and $u_2(z,t)$ of  Example \ref{exm:determ} are presented in Figure \ref{fig:test1_exact1} and Figure \ref{fig:test1_exact2} respectively. Computations have been carried out by MatLab$^{\copyright}$ R2019b \cite{Matlab} for Windows 10 home (64-bit) Intel(R) Core(TM) i5-8265u CPU, 1.60 GHz.

\begin{table}
	\centering
	\begin{tabular}{|c|ccc|}
	\hline 
	R 
	& RMSE$\left(u_1(z,t)\right)$ &  RMSE$\left(u_2(z,t)\right)$&  CPU, s\\
		\hline 
$ 5$	& $1.9717e-02$ & 	 $4.3759e-04$ & 	 $0.6406$ \\ 	
$ 10$	& $8.0645e-03$ & 	 $4.2821e-05$ & 	 $0.7188$ \\ 	
$ 15$	& $4.9310e-03$ & 	 $1.0860e-05$ & 	 $2.1094$ \\ 	
$ 20$	& $3.5510e-03$ & 	 $4.1179e-06$ & 	 $2.6406$ \\ 	
$ 25$	& $2.7871e-03$ & 	 $1.9559e-06$ & 	 $2.0938$ \\ 	
$ 30$	& $2.3023e-03$ & 	 $1.0771e-06$ & 	 $2.7969$ \\ 	 
\hline
	\end{tabular}
	\caption{Root mean square errors (RMSE) of numerical approximations of (\ref{eq:uR}) by the midpoint Riemann sum with fixed $h=0.05$ in the domain $(z,t)\in{[0,5]\times[0,1]}$ for the step-sizes $\Delta z=0.05$ and $\Delta t=0.01$.}
	\label{table:rmse_fixed_h}
\end{table}

\section{The Random Case: Numerical Solution and Numerical Convergence} \label{sec:random}

In previous Section \ref{sec:determ}, for a fixed realization corresponding to some $\xi_0 \in \Omega$, we have confirmed that Gaussian quadrature is not appropriate for approximating oscillatory integrals of the cosine inverse Fourier transform, see Example \ref{exm:determ}. Then, we experimented in Section \ref{sec:determ} the alternative of truncation combined with midpoint Riemann sum in the truncated domain $[0,R]$. In this section we only choose this technique to approximate the solution s.p. of the problem (\ref{eq:problem1})--(\ref{eq:problem4}), as well as the computation of the expectation and variance of the approximate s.p.

For the sake of clarity in the presentation, we recall some notation about s.p. used in previous paper \cite{Casaban2019a} and references therein.

Let $(\Omega, \mathfrak{F}, \mathbb{R})$ be a complete probability space, and let $L_p^{m \times n}(\Omega)$ be a set of all random  matrices $Y = (y_{ij})_{m \times n}$, whose entries $y_{ij}$ are r.v.'s satisfying

\begin{equation}
\norm{y_{ij}}_p = \left(\mathbb{E}[|y_{ij}|^p] \right)^{1/p} < + \infty, \quad p \geq 1,
\end{equation}
what means that $y_{ij} \in L_p(\Omega)$, where $\mathbb{E}[\cdot]$ denotes the expectation operator. The space of all random matrices of size $m \times n$, endowed with the matrix $p-$norm, $\left(L_p^{m \times n}(\Omega), \norm{\cdot}_p\right)$, defined by

\begin{equation}\label{eq:matrix_norm}
\norm{Y}_p = \sum_{i=1}^{m} \sum_{j=1}^{n} \norm{y_{ij}}_p, \quad \mathbb{E}[|y_{ij}|^p] < + \infty,
\end{equation}
is a Banach space. The definition of the matrix $p-$norm in (\ref{eq:matrix_norm}) can be extended to matrix s.p.'s $Y(t) = \left(y_{ij}(t) \right)_{m \times n} $ of $L_p^{m \times n}(\Omega)$, where now each entry $\left(y_{ij}(t) \right)$ is a s.p., that is, $y_{ij}(t)$ is a r.v. for each $t$. We say that the matrix s.p. $Y(t)$ lies in  $L_p^{m \times n}(\Omega)$, if $\left(y_{ij}(t) \right) \in L_p(\Omega)$ for every $1 \leq i \leq m$, $1 \leq j \leq n$.

The definitions of integrability, continuity, differentiability of a matrix function lying in $L_p^{m \times n}(\Omega)$ follows in a natural manner using matrix $p-$norm introduced in (\ref{eq:matrix_norm}). The case mean square corresponds to $p=2$, and mean four  $p=4$. One has $L_4^{m \times n}(\Omega) \subset L_2^{m \times n}(\Omega)$, see \cite{Casaban2019a}.

Consider a constant random matrix $L \in L_{2p}^{n \times n}(\Omega)$, $Y_0 \in L_{2p}^{n \times 1}(\Omega)$  and $C(s)$ lies in $L_{2p}^{n \times 1}(\Omega)$ and $2p-$integrable. Assume that random matrix $L=\left(\ell_{ij}\right)_{n \times n}$ satisfies the moment condition:

\begin{equation}\label{eq:moment_condition}
\mathbb{E}[|\ell_{ij}|^r] \leq m_{ij} (h_{ij})^r < + \infty, \quad \forall r \geq 0, \; \forall i,j: \; 1 \leq i,j \leq n.
\end{equation}

Then, by Section 3 of \cite{Casaban2019a}, the corresponding solution of the mean square initial value problem

\begin{equation}
Y'(s) = LY(s) + C(s), \quad Y(0) = Y_0, \; s \geq 0,
\end{equation}
is given by 
\begin{equation}
Y(s) = e^{Ls} \left( Y_0 + \int_{0}^{s}e^{-Lv} C(v)\ dv \right) .
\end{equation}

Now for the sake of clarity in the presentation we list the conditions of the random parabolic partial differential system (\ref{eq:problem1})--(\ref{eq:problem4}).  Firstly, random matrices

\begin{equation}
A(\xi) \in L_4^{2 \times 2}(\Omega), \quad B(\xi) \in L_8^{2 \times 2}(\Omega)\,,
\end{equation}
and satisfy the moment condition (\ref{eq:moment_condition}). In addition, assume that $g(t)(\xi) \in L_8^{2 \times 1}(\Omega)$ for each $t >0$ and it is continuous in $L_8^{2 \times 1}(\Omega)$. Furthermore assume that s.p. $f(z)(\xi)\in{L_8^{2 \times 1}(\Omega)}$ for each $z>0$ and it is absolutely integrable in  $L_8^{2 \times 1}(\Omega)$, then the random cosine Fourier transform of $f(z)(\xi)$, $F(\omega)(\xi)$, lies in $L_8^{2 \times 1}(\Omega)$. Assume the random spectral condition: There exists $b^{*} >0$ such that

\begin{equation}  \label{eq:randomSpectralCondition}
\inf_{\xi \in \Omega} \lambda_{\min} \left( \frac{B(\xi) + B (\xi)^T}{2}\right) \geq b^{*} >0.
\end{equation}

Then, the random initial value problem (see (\ref{eq:ode}))

\begin{equation}
\frac{d}{dt} V (t,\xi) (\omega) = (A(\xi)-\omega^2 B(\xi))V(t,\xi) (\omega) - B(\xi) g(t)(\xi), \quad  V(0,\xi)(\omega) = F(\omega),
\end{equation}
has the mean square solution

\begin{equation}  \label{eq:V}
V(t, \xi)(\omega) = e^{\left(A(\xi) - \omega^2 B(\xi)\right)\,t }
\left( F(\omega)(\xi) -\int_{0}^{t} e^{-\left(A(\xi) - \omega^2 B(\xi)\right)\, s}c(s)(\xi) ds\right), \quad c(s)(\xi) = B(\xi)g(s)(\xi),
\end{equation}
and $V(t,\xi)$ lies in $L_2^{2 \times 1} (\Omega)$. Using random cosine inverse transform, see (\ref{eq:u_def}), one has

\begin{equation}  \label{eq:SolExact}
u(z,t) (\xi) = \frac{2}{\pi} \int_{0}^{\infty} V(t,\xi)(\omega) \cos (\omega z) d \omega, \quad \xi \in \Omega.
\end{equation}

Truncation gives the approximate s.p.
\begin{equation}\label{eq:u_R_random}
u_R(z,t)(\xi) = \frac{2}{\pi}\int_{0}^{R} V(t,\xi)(\omega) \cos (\omega z) d \omega, \quad \xi \in \Omega.
\end{equation}

Using $N-$midpoint Riemann sum quadrature of $u_R(z,t)(\xi)$ one gets

\begin{equation}  \label{eq:Mid_u}
\text{Mid}_N\!\left(u_R(z,t) (\xi) \right) =\frac{2h}{\pi}  \sum_{j=0}^{N-1} V(t, \xi) \left(\omega_{j+ \frac{1}{2}}\right)\, \cos \left(\omega_{j+ \frac{1}{2}}z\right),  
\end{equation}
where $Nh = R$, $\omega_{j+ \frac{1}{2}} = \left(j+ \frac{1}{2} \right)h $.

Since the randomness of the integrand of (\ref{eq:u_R_random}) only affects to $V(t,\xi)(\omega)$, it is easy to check that 

\begin{equation}   \label{eq:EMid_N_R}
\mathbb{E} \left[ \text{Mid}_N \left(u_R(z,t) (\xi) \right) \right] =  \frac{2h}{\pi}  \sum_{j=0}^{N-1} \cos \left(\omega_{j+ \frac{1}{2}}z\right) 
\mathbb{E} \left[ V(t, \xi)\left(\omega_{j+ \frac{1}{2}}\right) \right]\,,
\end{equation}
and

\begin{equation} \label{eq:VarMid_N_R}
\begin{split}
\text{Var}\left[ \text{Mid}_N\! \left(u_R(z,t) (\xi) \right) \right] =
\mathbb{E} \left[ \left( \text{Mid}_N\! \left(u_R(z,t) (\xi) \right)\right) ^2\right] - \mathbb{E} \left[ \text{Mid}_N\! \left(u_R(z,t) (\xi) \right) \right]^2  
\\= \frac{4h^2}{\pi^2} \sum_{i=0}^{N-1} \sum_{j=0}^{N-1} \cos \left(\omega_{i+ \frac{1}{2}}z\right) \cos \left(\omega_{j+ \frac{1}{2}}z\right)  \text{cov} \left[ V(t,\xi) \left(\omega_{i+ \frac{1}{2}}\right), \,  V(t,\xi) \left(\omega_{j+ \frac{1}{2}}\right)\right], 
\end{split}
\end{equation}
where $\text{cov}[P,Q] = \mathbb{E}[PQ] -  \mathbb{E}[P]  \mathbb{E}[Q]$. 

Algorithm \ref{alg1} summarizes the steps to compute the approximations of the expectation and the standard deviation of the solution s.p. \eqref{eq:SolExact}

\begin{algorithm}
\caption{Procedure to compute the expectation and the standard deviation of the approximate solution s.p.  $\text{Mid}_N\left(u_{R}(z,t)(\xi)\right)$ \eqref{eq:Mid_u} of the problem \eqref{eq:problem1}--\eqref{eq:problem4}.}\label{alg1}
\begin{algorithmic}
\State Take random matrices $A(\xi)\in{L_4^{2 \times 2}(\Omega)}$ and $B(\xi)\in{L_8^{2 \times 2}(\Omega)}$ and check that all their entries satisfy the moment condition \eqref{eq:moment_condition}.
\State Take a continuous s.p. $g(t)(\xi)$ in $L_8^{2 \times 1}(\Omega)$ for $t>0$.
\State Take an absolutely integrable s.p. $f(z)(\xi)$ in $L_8^{2 \times 1}(\Omega)$ for $z>0$.
\State Check that the random spectral condition \eqref{eq:randomSpectralCondition} is verified.
\State Fix a point $(z,t)$, with $z>0$, $t>0$.
\State Choose the length of the truncation end-point $R$ and the number of subintervals $N$. 
\State Compute the step-size $h$ using the relationship $Nh=R$.
\For{$j=0$ to $j=N-1$}
\State compute the $N$-midpoints $\omega_{j+\frac{1}{2}}$ and the functions $\cos\left( \omega_{j+\frac{1}{2}} \,z \right)$.
\EndFor
\State Choose and carry out a number $K$ of realizations, $\xi=\{1, \ldots,K\}$, over the r.v.'s involve in the random matrices $A(\xi)$ and $B(\xi)$ and the s.p.'s $g(t)(\xi)$ and $f(z)(\xi)$. 
\For{each realization $\xi=1$ to $\xi=K$}
\State Compute the Fourier cosine transform of $f(z)(\xi)$: $F(\omega)(\xi)$.
\EndFor
\For{$j=0$ to $j=N-1$} 
\For{$\xi=1$ to $\xi=K$} 
\State Compute the deterministic expression \eqref{eq:V} with $\omega=\omega_{j+\frac{1}{2}}$. 
\EndFor
\State Compute the mean of the $K$ values obtained.
\EndFor
\State Compute the approximation of the expectation, $\mathbb{E}[\text{Mid}_N\left(u_{R}(z,t)(\xi)\right)]$ using expression \eqref{eq:EMid_N_R}.
\State Compute the approximation of the standard deviation, $\sqrt{\text{Var}[\text{Mid}_N\left(u_{R}(z,t)(\xi)\right)]}$ using expression \eqref{eq:VarMid_N_R}.
\end{algorithmic}
\end{algorithm}

In the next example we consider a random version of Example \ref{exm:determ} where the uncertainty in the computation of the exact values of coriolis and eddy viscosity is considered. Another uncertainty approach has been recently treated in \cite{Yosef2015}. 
\begin{exm}\label{exm:random}
We consider the random coupled parabolic problem (\ref{eq:problem1})--(\ref{eq:problem4}) with the initial a boundary conditions
\begin{equation} \label{eq:Boundary1Ex}
u(z,0) = f(z) = \begin{bmatrix}
	0\\0
	\end{bmatrix}, \quad \der{u}{z}(0,t) = \begin{bmatrix}
	-g(t)\\0
	\end{bmatrix}=\begin{bmatrix}
	-1
	\\0
	\end{bmatrix},
	\end{equation}
and the random matrix coefficients 
\begin{equation} \label{eq:Boundary2Ex}
	A(\xi) = \begin{bmatrix}
	0 & a(\xi) \\-a(\xi) & 0
	\end{bmatrix}, \quad 
	B(\xi)= \nu(\xi) I = \begin{bmatrix}
	\nu(\xi) & 0\\ 0 & \nu(\xi)
	\end{bmatrix},
	\end{equation}
where the r.v. $a(\xi)>0$ follows a Gaussian distribution of mean $\mu=2$ and standard deviation $\sigma=0.1$ truncated on the interval $[0.8, 1.2]$, that is $a(\xi) \sim N_{[0.8, 1.2]}{(2,0.1)}$, and the r.v. $\nu(\xi) >0$ has a gamma distribution of parameters $(4;2)$
truncated on the interval $[0.5, 1.5]$, that is $\nu(\xi) \sim Ga_{[0.5, 1.5]}{(4;2)}$. Both r.v.'s are considered independent ones.
Observe that the random matrices $A(\xi)\in{L_4^{2 \times 2}(\Omega)}$ and $B(\xi)\in{L_8^{2 \times 2}(\Omega)}$ and verifying the moment condition \eqref{eq:moment_condition} because the r.v.'s $a(\xi)$ and $b(\xi)$ are bounded. The function $g(t)=1$ involves in the boundary condition \eqref{eq:Boundary1Ex} is in $L_8^{2 \times 1}(\Omega)$ for each $t$ and is continuous too. Furthermore, the random spectral condition \eqref{eq:randomSpectralCondition} is satisfied because the r.v. $\nu(\xi)>0$. 
The exact solution of problem (\ref{eq:problem1})--(\ref{eq:problem4}) in its deterministic version is given by \eqref{eq:exact} with $g(s)=1$. Figure \ref{fig:MomentsExactSol} shows the numerical values of the expectation and the standard deviation of the exact solution s.p. \eqref{eq:problem1}--\eqref{eq:problem4}, \eqref{eq:Boundary1Ex}--\eqref{eq:Boundary2Ex} considering both at fixed time $t=1$. Computations have been carried out by Mathematica$^{\copyright}$ software version 11.3.0.0, \cite{Mathematica} for Windows 10Pro (64-bit) Intel(R) Core(TM) i7-7820X CPU, 3.60 GHz 8 kernels. In Table \ref{table:CPUsExactMoments} we show the timings (CPU time spent in the Wolfram Language kernel) to compute both statistical moments of the exact solution plotted in Figure \ref{fig:MomentsExactSol}.   

\begin{table}
	\centering
	\begin{tabular}{|cr|}
	\hline 
Statistical Moments &  CPU, s\\
	\hline 
$\mathbb{E}[u_1(z,1)]$             &  $18.2656$    \\
$\mathbb{E}[u_2(z,1)]$             &  $20.4531$    \\
$\sqrt{\mathrm{Var}[u_1(z,1)]}$    &  $7592.3800$    \\
$\sqrt{\mathrm{Var}[u_2(z,1)]}$    & $539.1560$     \\ 
\hline
	\end{tabular}
	\caption{Timings for computing the numerical values of the expectation and the standard deviation of the exact solution s.p. \eqref{eq:problem1}--\eqref{eq:problem4}, \eqref{eq:Boundary1Ex}--\eqref{eq:Boundary2Ex}
at  time $t=1$ in the spatial domain $0 \leq z \leq 5$ for $\Delta z=0.1$.}
	\label{table:CPUsExactMoments}
\end{table}

Numerical convergence of the statistical moments \eqref{eq:EMid_N_R}--\eqref{eq:VarMid_N_R} of the approximate solution s.p. \eqref{eq:Mid_u}
based on the midpoint Riemann sum quadrature and Monte Carlo technique is illustrated in the following way. Table \ref{table:rmse_Eu_fixed_R20N400} and Table \ref{table:rmse_DTu_fixed_R20N400} collect  the RMSEs for the numerical expectation and the standard deviation, respectively, for several number of realizations $\xi_i$ in the Monte Carlo method. In this first experiment  suitable values of truncation end-point $R$ and step-size $h$ suggested by the deterministic Example \ref{exm:determ} have been fixed. Figure \ref{fig:EyDTKsR20N400} illustrates the decreasing trend of the absolute errors of the approximations to the expectation, $\mathbb{E}[\text{Mid}_N\left(u_{R}(z,t)(\xi)\right)]$ \eqref{eq:EMid_N_R}, and the standard deviation, $\sqrt{\text{Var}[\text{Mid}_N\left(u_{R}(z,t)(\xi)\right)]}$ \eqref{eq:VarMid_N_R}, when the simulations $\xi_i$ by Monte Carlo  increase. If more precision is requiered it should be increased the values of parameters $R$ and $h$ (or $N$) rather than increasing the number of simulations $\xi_i$.

Secondly, by varying the length of $R$ for both $h$ and the number of realizations $\xi$ fixed, the computed RMSEs are shown in Table \ref{table:rmse_Eu_fixed_K1600h005Ni} and Table \ref{table:rmse_DTu_fixed_K1600h005Ni}. 
We have chosen $\xi=1600$ realizations due to the results obtained in the previous study are sufficiently accurate. Figure \ref{fig:ERsK} shows  how the approximate expectation, $\mathbb{E}[\text{Mid}_N\left(u_{R_i}(z,t)(\xi)\right)]$ \eqref{eq:EMid_N_R}, improves when the size of the integration domain $R_i$ increases.

The CPU times of the numerical experiments exhibit the efficience of the proposed method versus the long  times spent to compute the statistical moments  in Table \ref{table:CPUsExactMoments}.

\begin{table}
	\centering
	\begin{tabular}{|c|ccc|ccc|}
	\hline 
	$\xi_i$ 	& RMSE$\left(\mathbb{E}_{N}[u_{1,R}(z,1)(\xi_i)]\right)$  &  RMSE$\left(\mathbb{E}_{N}[u_{2,R}(z,1)(\xi_i)]\right)$           &  CPU, s\\
	\hline 
$200$   &  $4.9594e-03$  & $2.8818e-04$   & $0.1205$  \\
$400$   &  $5.6568e-03$  &  $1.2430e-03$  & $0.2500$ \\
$800$   &  $4.8149e-03$  &  $1.2268e-03$  &  $0.4531$ \\
$1600$  & $4.9603e-03$   &  $4.5712e-04$  &  $0.4844$\\
$3200$  & $5.0119e-03$   & $4.4230e-04$   &  $1.1406$ \\
$6400$  & $4.7532e-03$   &  $8.8098e-04$  &  $1.2188$ \\
$12800$ & $4.8759e-03$   & $1.2562e-04$   &  $5.2656$ \\		\hline
	\end{tabular}
	\caption{Root mean square errors (RMSEs) of the numerical approximations of the expectation  \eqref{eq:EMid_N_R} for the  the solution s.p. \eqref{eq:SolExact}  at $t=1$ in  $0 \leq z \leq 5$  with $\Delta z=0.1$ varying the number of realizations $\xi_i$ for $R = 20$ and $h=0.05$ ($N=400$). }
	\label{table:rmse_Eu_fixed_R20N400}
\end{table}

\begin{table}
	\centering
	\begin{tabular}{|c|ccc|ccc|}
	\hline 
	$\xi_i$ 	& RMSE$\left(\sqrt{\text{Var}_{N}[u_{1,R}(z,1)(\xi_i)]}\right)$ &  RMSE$\left(\sqrt{\text{Var}_{N}[u_{2,R}(z,1)(\xi_i)]}\right)$          &  CPU, s\\
	\hline 
$200$   & $1.8818e-03$       & $1.3906e-03$  & $0.6094$ \\
$400$   & $2.7378e-04$     & $1.9727e-04$ &  $0.9531$  \\
$800$   & $3.8938e-03$     & $2.3218e-03$  & $1.9844$  \\
$1600$  & $5.0547e-04$     & $3.2676e-04$ & $6.9531$ \\
$3200$  &  $3.1779e-04$     & $1.1304e-04$  & $19.0313$  \\
$6400$  & $6.6312e-04$       & $3.6929e-04$  & $55.1094$   \\
$12800$ & $2.9489e-04$      & $1.5663e-04$ & $186.484$   \\		\hline
	\end{tabular}
	\caption{Root mean square errors (RMSE) of the numerical approximations of  the standard deviation \eqref{eq:VarMid_N_R} for the solution s.p. \eqref{eq:SolExact}  at $t=1$ in $0 \leq z \leq 5$ with $\Delta z=0.1$ varying the number of realizations $\xi_i$ for $R = 20$ and $h=0.05$ ($N=400$).  }
	\label{table:rmse_DTu_fixed_R20N400}
\end{table}

\begin{table}
	\centering
	\begin{tabular}{|c|ccc|ccc|}
	\hline 
	$R_i$ 	& RMSE$\left(\mathbb{E}_{N_i}[u_{1,R_i}(z,1)(\xi)]\right)$  &  RMSE$\left(\mathbb{E}_{N}[u_{2,R_i}(z,1)(\xi)]\right)$          &  CPU, s\\
	\hline 
$5$   &  $2.2870e-02$   & $7.0643e-04$  &  $0.1719$ \\
$10$  &  $1.0113e-02$   & $4.6269e-04$  &  $0.3281$\\
$15$  &  $6.6053e-03$   &  $4.5785e-04$ &  $0.5938$ \\
$20$  &  $4.9603e-03$   &  $4.5712e-04$ &  $0.6406$\\
$25$  &  $3.9572e-03$   &  $4.5689e-04$ &  $1.3750$ \\
	\hline
	\end{tabular}
	\caption{Root mean square errors (RMSE) of the numerical approximations of the expectation  \eqref{eq:EMid_N_R} for the solution s.p. \eqref{eq:SolExact}  at $t=1$ in  $0 \leq z \leq 5$ with  $\Delta z=0.1$. The number of realizations $\xi=1600$ and $h=0.05$ are fixed and the size of the integration domain $R_i$ varies. }
	\label{table:rmse_Eu_fixed_K1600h005Ni}
\end{table}

\begin{table}
	\centering
	\begin{tabular}{|c|ccc|ccc|}
	\hline 
	$R_i$ 	& RMSE$\left(\sqrt{\text{Var}_{N}[u_{1,R}(z,1)(\xi)]}\right)$   &  RMSE$\left(\sqrt{\text{Var}_{N}[u_{2,R_i}(z,1)(\xi)]}\right)$         &  CPU, s\\
	\hline 
$5$    & $5.0618e-04$     & $3.4636e-04$  & $3.3594$ \\
$10$  & $5.0549e-04$     & $3.2624e-04$  &  $13.1406$ \\
$15$  &  $5.0547e-04$    & $3.2662e-04$  &  $34.6875$ \\
$20$  & $5.0547e-04$     & $3.2676e-04$  &  $80.4688$\\
$25$  & $5.0547e-04$    & $3.2681e-04$  &   $123.1250$\\
	\hline
	\end{tabular}
	\caption{Root mean square errors (RMSEs) of the numerical approximations of the expectation  standard deviation \eqref{eq:VarMid_N_R} for the solution s.p.\eqref{eq:SolExact}  at $t=1$ in  $0 \leq z \leq 5$ with  $\Delta z=0.1$. The number of realizations $\xi=1600$ and $h=0.05$ are fixed and  the size of the integration domain $R_i$ varies. }
	\label{table:rmse_DTu_fixed_K1600h005Ni}
\end{table}


\begin{figure}[h]
\centering
  \includegraphics[width=1\textwidth]{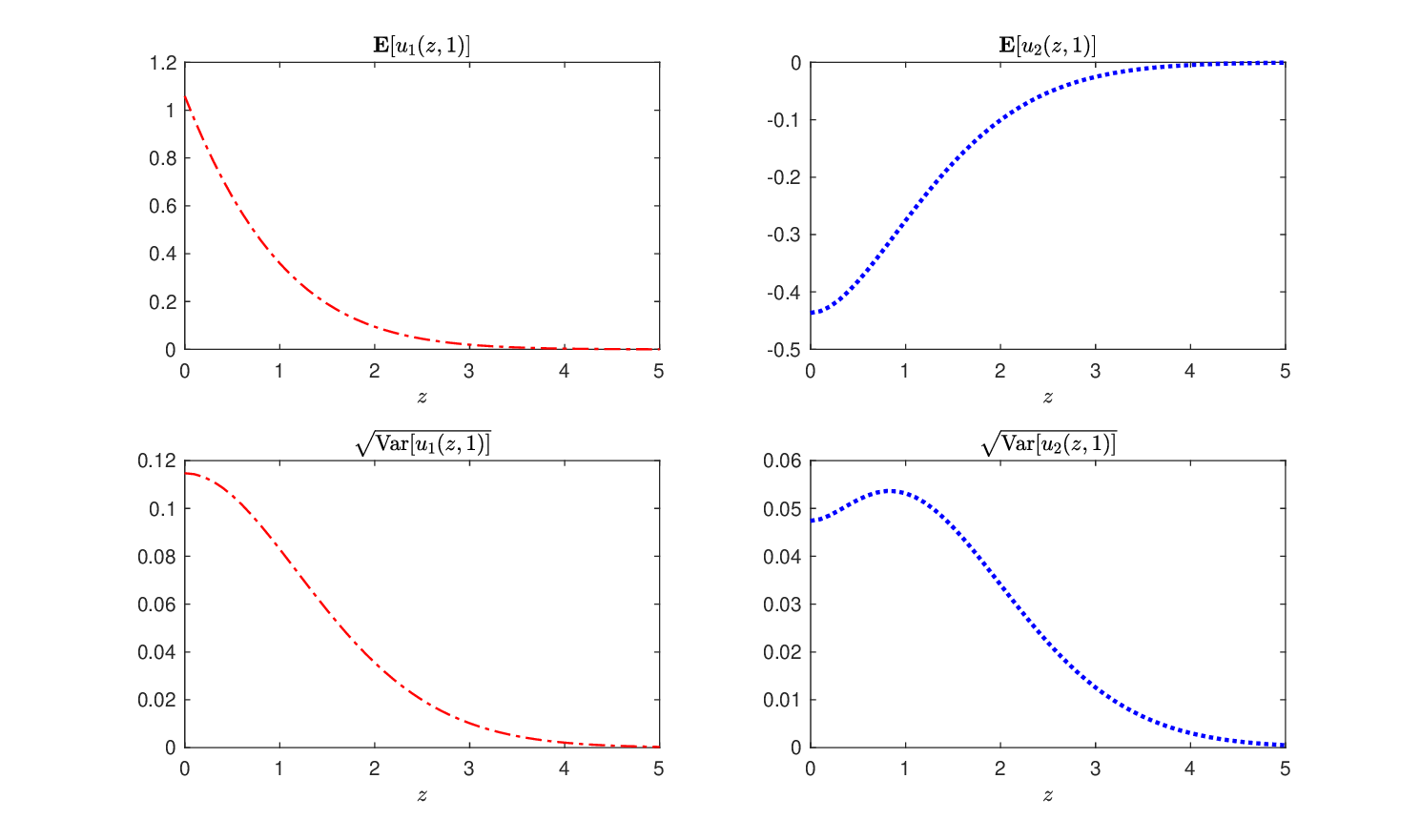}
 \caption{\label{fig:MomentsExactSol} Expectation and the standard deviation at the time instant $t=1$ of the exact solution s.p, $[u_1(z,t),u_2(z,t)]$, for the random coupled parabolic problem (\ref{eq:problem1})--(\ref{eq:problem4}), \eqref{eq:Boundary1Ex}--\eqref{eq:Boundary2Ex}, considering the r.v.'s $a(\xi) \sim N_{[0.8, 1.2]}{(2,0.1)}$ and $\nu(\xi) \sim Ga_{[0.5, 1.5]}{(4;2)}$, and the spatial domain $z\in[0,5]$ with step size $\Delta z=0.1$.}
\end{figure}

\begin{figure}[h]
\centering
  \includegraphics[width=1\textwidth]{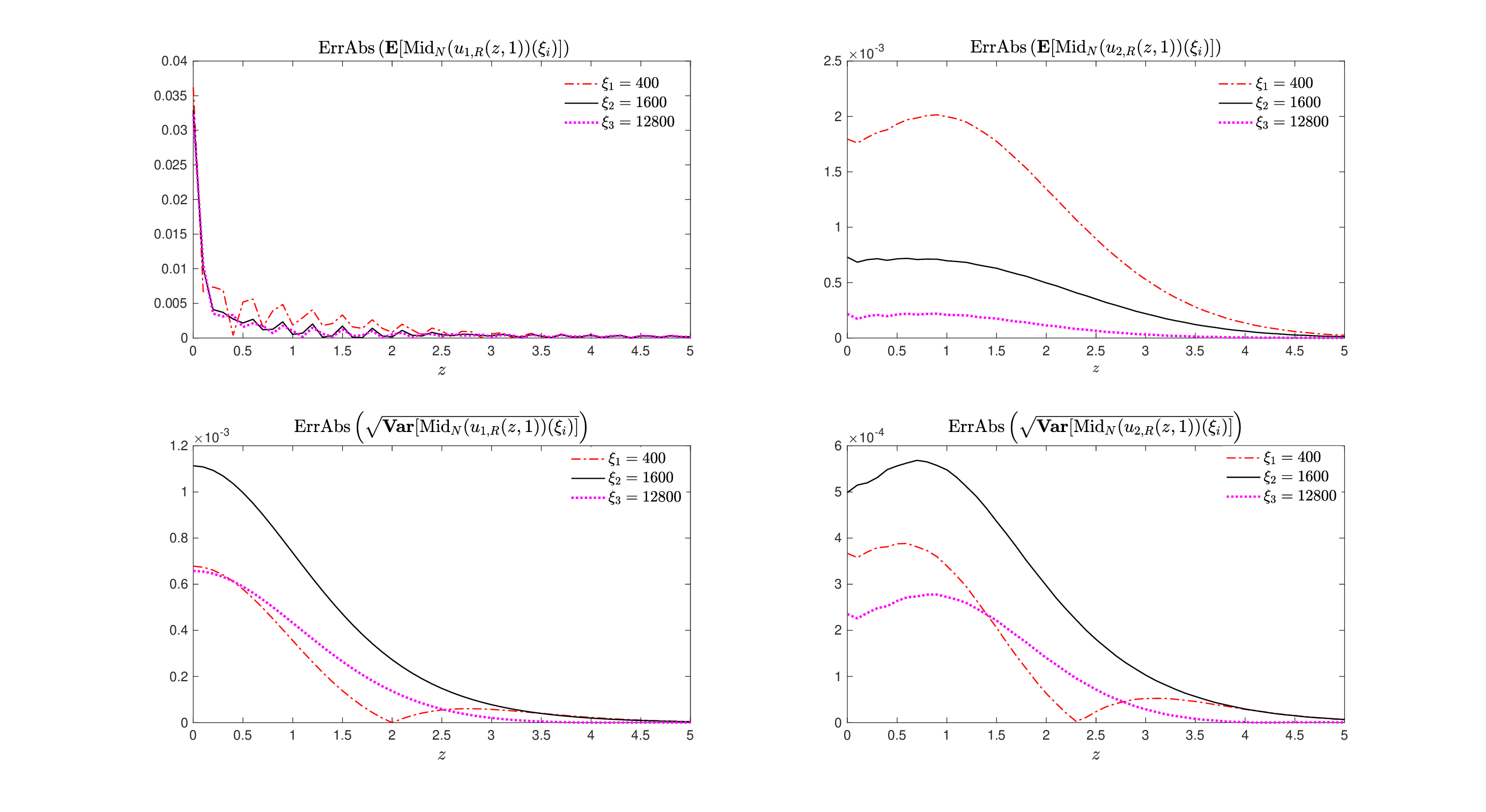}
 \caption{\label{fig:EyDTKsR20N400} Absolute errors of the expectation and the standard deviation for both components of the approximate solution s.p. \eqref{eq:Mid_u} at $t=1$ fixing $R=20$ and $h=0.05$ ($N=400$) in \eqref{eq:EMid_N_R}--\eqref{eq:VarMid_N_R} but varing the number of simulations $\xi_i=\{ 400, 1600, 12800\}$. The spatial domain is $z\in[0,5]$ with $\Delta z =0.1$.}
\end{figure}

\begin{figure}[h]
\centering
  \includegraphics[width=1\textwidth]{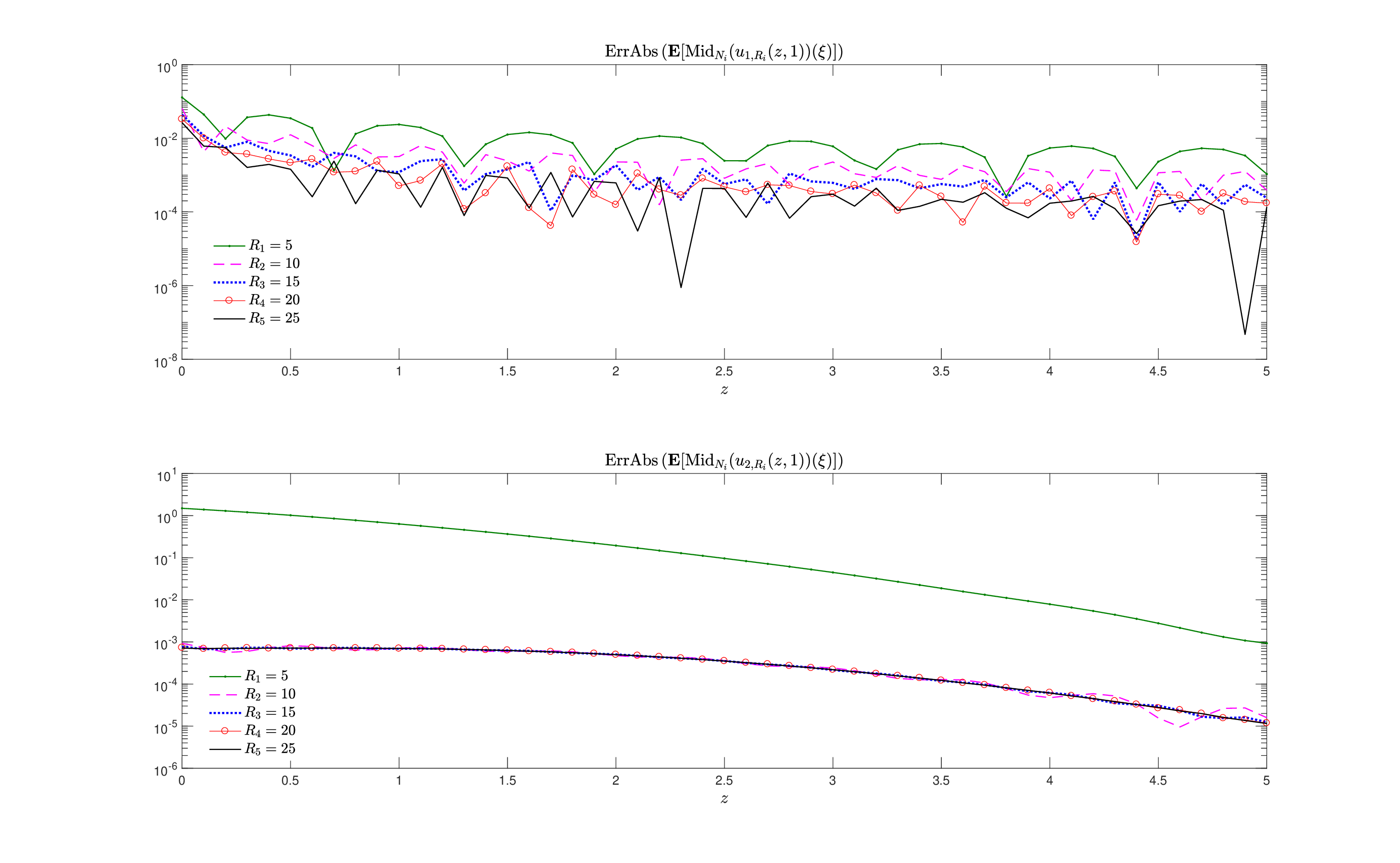}
 \caption{\label{fig:ERsK} Absolute errors of the expectation and the standard deviation for both components of the approximate solution s.p. \eqref{eq:Mid_u} at $t=1$ with logarithmic scale for the y-axis. The number of simulations $\xi=1600$ and the step-size $h=0.05$ in \eqref{eq:EMid_N_R}--\eqref{eq:VarMid_N_R} are fixed but the size of the integration domain, $R$, varies $R_i=\{5,10,15,20,25\}.$ The spatial domain is $z\in[0,5]$ with $\Delta z =0.1$ and $N_{i}=\{100, 200,300,400,500 \}$ so that $R=Nh$.}
\end{figure}

\end{exm}

\section{Conclusions}
This paper shows that the integral transform method combined with numerical integration and Monte Carlo technique is useful to deal with random systems of partial differential equations. Although here we consider parabolic type systems and cosine Fourier transform, the ideas are applicable to other systems and other integral transforms. The numerical integration must consider the possible oscillatory nature of the involved integrals. This approach is a manageable alternative to deal the computational complexity derived from the treatment of random models versus the iterative methods.

\section*{Acknowledgements}
This work has been partially supported by the Ministerio de Ciencia, Innovación y Universidades, Spanish grant MTM2017-89664-P.

\bibliography{bib-parabolic}

\end{document}